\documentclass[Leqno,12pt]{amsart}
\usepackage{amsmath}
\usepackage{amsxtra}
\usepackage{amscd}
\usepackage{amsthm}
\usepackage{amsfonts}
\usepackage{amssymb}
\usepackage{eucal}
\usepackage{mathrsfs}
\usepackage{color, graphics}
\usepackage[all]{xy}
\usepackage{hyperref}
\usepackage{pstricks,pst-plot}

\def\C{{\mathbb C}}

\def\E{{\mathbb E}}

\def\P{{\mathbb P}}

\def\Z{{\mathbb Z}}

\title[Effects of property $N_p$]{Some effects of property $N_p$ on the higher normality and
defining equations of nonlinearly normal varieties}

\author[S. Kwak]{Sijong Kwak$^{\star}$}
\address%[S. Kwak]
{Sijong Kwak : Department of Mathematics, Korea Advanced Institute
of Science and Technology, 373-1 Gusung-dong, Yusung-gu, Taejon,
Korea} \email{sjkwak@math.kaist.ac.kr}
\thanks{$\star$ Supported by grant No. R01-2002-000-00151-0
from the Research Program of  Korea Science \&  Engineering
Foundation(KOSEF)}

\author[E. Park]{Euisung Park$^{\star \star}$}
\address%[E. Park]
{Euisung Park : Department of Mathematics, Korea Advanced
Institute of Science and Technology, 373-1 Gusung-dong, Yusung-gu,
Taejon, Korea} \email{puspus@kaist.ac.kr}
\thanks{$\star \star$ Supported by Korea Research Foundation Grant (KRF-2002-070-C00003)}

\begin{document}

\newtheorem{theorem}{Theorem}[section]
\newtheorem{definition}{Definition}[section]
\newtheorem{lemma}[theorem]{Lemma}
\newtheorem{remark}{Remark}
\newtheorem{con}{Conjecture}
\newtheorem{ex}{Example}
\newtheorem{proposition}[theorem]{Proposition}
\newtheorem{examples}[theorem]{Examples}
\newtheorem{corollary}[theorem]{Corollary}

\thispagestyle{empty} \maketitle

\begin{abstract}
For a smooth projective variety $X\subset \P^r$ embedded by the
complete linear system, Property $N_p$ has been studied for a long
time(\cite{Butler}, \cite{Green}, \cite{GL}, \cite{EL} etc). On
the other hand, Castelnuovo-Mumford regularity conjecture and
related problems have been focused for a projective variety which
is not necessarily linearly normal(\cite{BM}, \cite{GLP},
 \cite{K}, \cite{L}, \cite{Noma} etc). This paper aims to explain
the influence of Property $N_p$ on higher normality and defining
equations of a smooth variety embedded by sub-linear system. Also
we prove a claim about Property $N_p$ of surface scrolls which is
a generalization of Green's work in \cite{Green} about Property
$N_p$ for curves.
\end{abstract}

\tableofcontents
 \setcounter{page}{1}

\section{Introduction}
\noindent When a projective variety $X$ is embedded in a
projective space $\P^r$, there are various natural interesting
questions regarding the syzygy modules of the saturated ideal
$I_X$ and the finitely generated graded module $\oplus_{\ell \in
\Z} H^0 (X, \mathcal{O}_X (\ell))$. For a linearly normal variety
$X \subset \P^r$, finding conditions to guarantee that $X
\hookrightarrow \P H^0 (X,L)$ is projectively normal or cut out by
quadrics has expanded to the condition for numerical type of
higher syzygies of the homogeneous coordinate ring $S(X)$, i.e.
Property $N_p$ (Castelnuovo\cite{C}, Mumford\cite{M2},
Green\cite{Green}, Ein-Lazarsfeld\cite{EL}, etc). On the other
hand, for a projective variety $X \subset \P^r$ which is not
necessarily linearly normal, the classical problem of finding
upper bound of $n$ satisfying
\begin{enumerate}
\item[(1)]
Hypersurfaces of degree $n$ cut out a complete linear system on
$X$.
\item[(2)]
$X$ is cut out in $\P^r$ by hypersurfaces of degree $n$, and the
homogeneous ideal of $X$ is generated in degrees $\geq n$ by its
component of degree $n$.
\end{enumerate}
has been reformulated as finding the upper bound of the regularity
of the ideal sheaf $\mathcal{I}_X$ of $X\subset \P^r$ in terms of
degree, codimension and other projective invariants
(Castelnuovo\cite{C}, Mumford\cite{M}, Eisenbud-Goto\cite{EG},
etc).

For precise statements, we give notations and recall definitions.
Let $X$ be a complex projective variety of dimension $n$ and let
$\mathcal{L}$ be a very ample line bundle on $X$. Consider an
embedding $X \subset \P(V)$ where $V \subset H^0 (X,\mathcal{L})~~
$ is a subsystem of dimension $r+1$. We use the following usual
notations:
\begin{enumerate}
\item[$\bullet$] $S=Sym ^{\bullet} V$ : the homogeneous coordinate ring of
$\P(V)=\P^r$
\item[$\bullet$] $I_X \subset S$ : the homogeneous ideal of $X$
\item[$\bullet$] $S(X) = S/I_X$ : the homogeneous coordinate ring of
$X$
\item[$\bullet$] $\mathcal{O}_{\P^r}, \mathcal{O}_X$ : the structure sheaves of $\P^r$ and
$X$, respectively
\item[$\bullet$] $\mathcal{I}_X \subset \mathcal{O}_{\P^r}$ : the sheaf
of ideals of $X$
\end{enumerate}
Now let us consider the minimal free resolution of $S(X)$:
\begin{equation*}
0 \rightarrow L_{r+1} \rightarrow \cdots \rightarrow L_i
\stackrel{\varphi_i}{\rightarrow} L_{i-1} \rightarrow \cdots
\rightarrow L_1 \stackrel{\varphi_1}\rightarrow L_0
\stackrel{\varphi_0}{\rightarrow} S(X) \rightarrow 0
\end{equation*}
where $L_i$, as free graded $S$-module, can be written as
\begin{equation*}
L_i = \oplus_j S^{k_{i,j}}(-i-j).
\end{equation*}
Note that $L_0 =S$ and the image of $\varphi_1$ is $I_X$. First
recall the definitions of Property $N_p$ and Castelnuovo-Mumford
regularity:

\begin{definition}  $(1)$ For positive integer $p\ge 1$,
$\mathcal{L}$ is said to satisfy Property $N_p$ if for $V=H^0
(X,\mathcal{L})$, a linearly normal embedding $X \subset \P(V)$ is
projectively normal, i.e., $S(X)= \oplus_{\ell \in \Z} H^0
(X,\mathcal{L}^{\ell})$, and
\begin{equation*}
k_{i,j} = 0~~~~\mbox{for all}~~~~0\leq i \leq
p~~~~\mbox{and}~~~~j\geq 2.
\end{equation*}
$(2)$ $X \subset \P(V)$, is said to be $m$-regular if we have
$k_{i,j} =0$ for all $j\ge m$, that is, one can find a basis of
$L_i$ with elements of degree at most $m+i-1$.
\end{definition}

Property $N_p$ means that the first $(p+1)$-th modules of syzygies
of $S(X)$ are as simple as possible or more precisely the
resolution is linear until the $p$-th stage. So, remark that
$\mathcal{L}$ satisfies Property $N_0$ if and only if
$\mathcal{L}$ is normally generated(i.e., projectively normal),
$\mathcal{L}$ satisfies Property $N_1$ if and only if
$\mathcal{L}$ satisfies Property $N_0$ and the homogeneous ideal
is generated by quadrics, $\mathcal{L}$ satisfies Property $N_2$
if and only if $\mathcal{L}$ satisfies Property $N_1$ and the
relations among the quadrics are generated by the linear relations
and so on.

 For Castelnuovo-Mumford regularity, there is another definition using vanishing
of certain cohomology groups due to D. Mumford. In fact, as D.
Mumford defined, $X$ is said to be $m$-regular if $H^i (\P^r ,
\mathcal{I}_X (m-i))=0$ for every $i\ge 1$. Recall (Lecture 14,
\cite{M}) that if $X$ is $m$-regular, then $X \subset \P^r$ is
$k$-normal for all $k\geq m-1$ and $\mathcal{I}_X (m)$ is globally
generated. Later D. Eisenbud and S. Goto gave the definition of
regularity for graded $S$-modules as above and proved the
equivalence of them \cite{EG}.\\

Now we generalize Property $N_p$ to non-complete embedding. Let
$X$ be a complex projective variety of dimension $n$ and let
$\mathcal{L}$ be a very ample line bundle on $X$. For an
isomorphic embedding
\begin{equation*}
X \hookrightarrow \P(V)
\end{equation*}
where $V \subset H^0 (X,\mathcal{L})$ a subsystem of codimension
$t$, let $S$ be the homogeneous coordinate ring of $\P(V)$. In
order to adjust the Koszul cohomology methods to noncomplete
embeddings, let us consider the graded $S$-module $E =
\oplus_{\ell \in \Z} H^0 (X, \mathcal{L}^{\ell})$. We consider the
situations in which the first few modules of syzygies of $E$ are
as simple as possible:

\begin{definition} Assume that $E$ has a minimal free resolution
of the form
\begin{equation*}
\cdots \stackrel{\varphi_{i+1}}{\rightarrow} \oplus_j S
^{k_{i,j}}(-i-j) \stackrel{\varphi_i}{\rightarrow}  \cdots
\stackrel{\varphi_1}{\rightarrow} \oplus_j S^{k_{0,j}}(-j)
\stackrel{\varphi_0}{\rightarrow} E \rightarrow 0.
\end{equation*}
For an integer $p \geq 0$, the embedding $X  \hookrightarrow
\P(V)$ is said to satisfy $N^{S}_p$ property if $k_{i,j} = 0$ for
$0 \leq i \leq p$ and  $j \geq 2$.
\end{definition}

\noindent Property $N^S _p$ means that $E$ admits a minimal free
resolution of the form
\begin{eqnarray*}
\cdots \rightarrow S^{m_p}(-p-1) \rightarrow \cdots \rightarrow
S^{m_2}(-3) \rightarrow  S^{m_1}(-2) \rightarrow S \oplus S^t (-1)
\rightarrow E \rightarrow 0.
\end{eqnarray*}
See also the following table of Betti
numbers:\\
\begin{table}[hbt]
\begin{center}
\begin{tabular}{|c|c|c|c|}\hline
                   &$k_{0,0}=1$     & $k_{0,1}$    &$k_{i,j}=0$        \\\cline{2-3}
                   &                & $k_{1,1}$      & for all          \\\cline{3-3}
    $k_{i,j}=0$    &   $k_{i,0}=0$  & $\vdots$       & $0 \leq i \leq p$ \\\cline{3-3}
for all $j\leq -1$ &for all $i\geq1$& $k_{p,1}$      & and
$j\geq2$
\\\cline{3-4}
                   &                & $\vdots$       &$\vdots$ \\\hline
\end{tabular}
\end{center}
\caption{Betti numbers when Property $N^{S}_p$ holds.}
\end{table}\\
Here $k_{0,1} =t$ and hence for linearly normal embedding,
Property $N^{S}_p$ is equal to Property $N_p$. This paper is
intended to investigate the following three problems:\\
\begin{enumerate}
\item[(1)] the influence of Property $N_p$ on higher syzygies of $X$ embedded by
subsystems( i.e., Property $N^{S}_q$ for $q \leq p-1$)
\item[(2)] the influence of Property $N^{S}_1$ on higher normality and defining
equations
\item[(3)] Property $N_p$ for surface scrolls (which is a generalization of
Green's ``$2g+1+p$" theorem for curves)\\
\end{enumerate}
Precisely, we give answers for these problems in Theorem
\ref{thm:regularity}, Theorem \ref{thm:effect} and Theorem
\ref{thm:surfacescroll}, respectively.

\begin{theorem}\label{thm:regularity}
Let $X$ be a smooth complex projective variety and let
$\mathcal{L} \in \mbox{Pic}X$ be a very ample line bundle. For an
embedding $ X \hookrightarrow \P(V)$ given by a subsystem $V
\subset H^0 (X,\mathcal{L})~~ $ of codimension $t$, if it
satisfies Property $N^{S}_1$, then $X \hookrightarrow \P(V)$ is
\begin{enumerate}
\item[$(1)$] $k$-normal for all $k \geq t+1$,
\item[$(2)$] max$\{m+1,t+2\}$-regular and,
\item[$(3)$] cut out by hypersurfaces of degree $t+2$, and the
homogeneous ideal of $X$ is generated in degrees $\geq t+2$ by its
component of degree $t+2$
\end{enumerate}
where $m$ is the regularity of $\mathcal{O}_X$ with respect to
$\mathcal{L}$.
\end{theorem}

\begin{theorem}\label{thm:effect}
Let $X$ be a smooth complex projective variety and let
$\mathcal{L} \in \mbox{Pic}X$ be a very ample line bundle such
that
\begin{equation*}
H^1 (X,\mathcal{L}^{j})=0~~~~ \mbox{for all}~~~~j\geq2.
\end{equation*}
If $ \mathcal{L}$ satisfies Property $N_p$, then for every
embedding $X \hookrightarrow \P(V )$ given by a subsystem $V
\subset H^0 (X,\mathcal{L})~~ $ of codimension $0 \leq t \leq p$,
Property $N^{S}_{p-t}$ holds. In particular, for $0 \leq t \leq
p-1$, $X \hookrightarrow \P(V)$ is
\begin{enumerate}
\item[$(1)$] $k$-normal for all $k \geq t+1$,
\item[$(2)$] max$\{m+1,t+2\}$-regular and,
\item[$(3)$] cut out by hypersurfaces of degree $t+2$, and the
homogeneous ideal of $X \hookrightarrow \P(V)$ is generated in
degrees $\geq t+2$ by its component of degree $t+2$
\end{enumerate}
where $m$ is the regularity of $\mathcal{O}_X$ with respect to
$\mathcal{L}$.
\end{theorem}

\noindent In brief, our main Theorems \ref{thm:regularity} and
\ref{thm:effect} imply that Property $N_p$ for some $p \geq 1$
guarantees a good behavior of higher normality,
Castelnuovo-Mumford regularity and optimal upper bound of the
degree of polynomials generating $I_X$ for an embedding by the
subsystem of codimension less than $p$. Furthermore we can apply
our main theorems to almost all known smooth projective varieties
with Property $N_p$ because the cohomological condition in Theorem
\ref{thm:effect} ($H^1 (X,\mathcal{L}^\ell)=0$ for $\ell \geq 2$)
is satisfied for all those cases (For details, see $\S 5$). We
also remark that Theorem \ref{thm:effect} can be applied to almost
all special curves since we do not assume $H^1 (X,\mathcal{L})=0$.

To get the upper bound of degree of minimal generators for the
homogeneous ideal, it has been canonical to investigate
Castelnuovo-Mumford regularity of the ideal sheaf. Along this
point of view, Theorem \ref{thm:regularity} implies surprisingly
that the degree bound can be much smaller than the optimal
regularity bound. In particular, this degree bound is also optimal
since there exists $(t+2)$-secant line to $X$ obtained from linear
projections with special centers. Also the bound of higher
normality and regularity can not be refined since in the case of
ruled scrolls over curves, the higher normality is equivalent to
regularity. In general, regularity or normality is worse as $t$ is
getting larger which is geometrically related to the existence of
a higher multisecant line. Remark that $X$ fails to be $k$-regular
if it admits $(k+1)$ secant line.\\

\noindent {\bf Remark 1.} In \cite{Birkenhake}, Christina
Birkenhake generalize Property $N_p$ in a different way. Let
$\mathcal{L}$ be a very ample line bundle on a smooth projective
variety $X$ and $V \subset H^0 (X,\mathcal{L})$ a subsystem
defining an embedding $X \hookrightarrow \P (V)$. For the
homogeneous coordinate ring $S = Sym^{\bullet} V$ of $\P (V)$,
define the graded $S$-module
\begin{equation*}
R := k \oplus V \oplus H^0 (X,\mathcal{L}^2) \oplus H^0
(X,\mathcal{L}^3) \oplus \cdots.
\end{equation*}
Then $X \hookrightarrow \P (V)$ is said to satisfy Property
$\tilde{N_p}$ if the minimal free resolution of $R$ is linear
until the $p$-th stage. That is, $R$ admits a minimal free
resolution of the form
\begin{eqnarray*}
\cdots \rightarrow S^{m_p}(-p-1) \rightarrow \cdots \rightarrow
S^{m_2}(-3) \rightarrow  S^{m_1}(-2) \rightarrow S \rightarrow R
\rightarrow 0.
\end{eqnarray*}
Therefore Property $\tilde{N_0}$ holds if and only if
$k$-normality holds for all $k \geq 2$, Property $\tilde{N_1}$
holds if and only if Property $\tilde{N_0}$ holds and the
homogenous ideal is generated by quadrics, Property $\tilde{N_2}$
if and only if Property $\tilde{N_1}$ holds and the relations
among the quadrics are generated by the linear relations and so
on. Clearly $S(X) \subseteq R \subseteq E$ and in particular $R =
E$ if and only if $V = H^0 (X,L)$. However, Koszul cohomological
method developed in \cite{Green} does not work well in the case of
Property $\tilde{N_p}$ as Alberto Alzati and Francesco Russo
mentioned in their paper \cite{AR}. To our knowledge, this is
because the graded module $R$ is not saturated. Indeed, for a very
ample line bundle $\mathcal{L}$ such that $H^1
(X,\mathcal{L}^j)=0$ for all $j \geq 2$ and Property $N_2$ holds,
every isomorphic one point projection satisfies Property $N^S _1$
and $k$-normality for all $k \geq 2$. But, there exists one point
projection admitting trisecant line in many cases and hence its
image cannot be cut out by quadrics ideal-theoretically, which
also supports the argument in \cite[Example 4.4]{AR}.
\\

\noindent {\bf Remark 2.} In \cite{AR}, Alberto Alzati and
Francesco Russo gave a necessary and sufficient condition that the
isomorphic projection of a $k$-normal variety to remain
$k$-normal. As an application, they also proved that if
$\mathcal{L}$ satisfies Property $N_2$, then every isomorphic one
point projection satisfies $k$-normality for all $k \geq 2$
(without the assumption that $H^1 (X,\mathcal{L}^j)=0$ for
all $j \geq 2$). For details, see Theorem 3.2 and Corollary 3.3 in \cite{AR}.\\

In $\S 4$ we study higher syzygies of linear series on surface
scrolls. Indeed we prove a criterion for Property $N_p$ of surface
scrolls and then apply Theorem \ref{thm:effect} to investigate the
geometry of subsystems of very ample line bundles giving scroll
embedding. We will follow the notation and terminology of R.
Hartshorne's book \cite{H}, V $\S 2$. Let $C$ be a smooth
projective curve of genus $g$ and let $\mathcal{E}$ be a vector
bundle of rank $2$ on $C$ which is normalized, i.e., $H^0
(C,\mathcal{E}) \neq 0$ while $H^0 (C,\mathcal{E} \otimes
\mathcal{O}_C (D))=0$ for every divisor $D$ of negative degree. We
set
\begin{equation*}
\mathfrak{e}=\wedge^2 \mathcal{E}~~~~\mbox{and}~~~~e = -
\mbox{deg}(\mathfrak{e}).
\end{equation*}
Let $X = \P_C (\mathcal{E})$ be the associated ruled surface with
projection morphism $\pi : X \rightarrow C$. We fix a minimal
section $C_0$ such that $\mathcal{O}_X (C_0)=\mathcal{O}_{\P_C
(\mathcal{E})} (1)$. For $\mathfrak{b} \in \mbox{Pic}C$,
$\mathfrak{b}f$ denote the pullback of $\mathfrak{b}$ by $\pi$.
Thus any element of $\mbox{Pic}X$ can be written
$aC_0+\mathfrak{b}f$ with $a\in \Z$ and $\mathfrak{b} \in
\mbox{Pic}C$ and any element of $\mbox{Num}X$ can be written $aC_0
+bf$ with $a,b \in \Z$ and scroll means that $a=1$. Precisely our
result is the following:

\begin{theorem}\label{thm:surfacescroll}
Let $\mathcal{L} =C_0 +\mathfrak{b}f \in \mbox{Pic}(X)$ be a line
bundle in the numerical class of $C_0 +bf$ such that
\begin{equation*}
b = \mbox{max} \{ 3g,3g+ e \} + p \quad \mbox{for some} \quad p
\geq 0.
\end{equation*}
Then $\mathcal{L}$ has Property $N_p$.
\end{theorem}

\noindent It is easily checked that $\mathcal{O}_X$ is $2$-regular
with respect to $\mathcal{L}$. So by applying Theorem
\ref{thm:effect} to this result, we obtain the following:

\begin{corollary}
Under the same situation as in Theorem \ref{thm:surfacescroll},
consider the embedding $X \hookrightarrow \P(V)$ given by a
subsystem $V \subset H^0 (X,\mathcal{L})~~ $ of codimension $t$.
If $0 \leq t \leq p-1$, then
$X \hookrightarrow \P(V) ~~~~ \mbox{is} \begin{cases} \mbox{$k$-normal for all $k \geq t+1$ and,}\\
\mbox{ max$\{3,t+2\}$-regular}.
\end{cases}$\\
In particular, $X \hookrightarrow \P(V)$ is cut out by
hypersurfaces of degree $t+2$, and the homogeneous ideal of $X$ is
generated in degrees $\geq t+2$ by its component of degree $t+2$.
\end{corollary}

The organization of this paper is as follows. In $\S 3$, we
investigate the geometric meaning of syzygies of the module $E$.
Precisely we prove Theorem \ref{thm:regularity} and Theorem
\ref{thm:effect}. In $\S 4$, we give a concrete statement and
proof of Property $N_p$ for surface scrolls. $\S 5$ is devoted to
apply Theorem \ref{thm:regularity} and Theorem \ref{thm:effect}.
In particular we reprove some old results about regularity and
Property $N_p$.\\

\section{Notations and Conventions}
Throughout this paper the following is assumed.\\\\
$(1)$ All varieties are defined over the complex number field $\C$.\\\\
$(2)$ For a finite dimensional $\C$-vector space $V$, $\P(V)$ is
the projective space of one-dimensional quotients of $V$.\\\\
$(3)$ When a variety $X$ is embedded in a projective space, we
always assume that it is non-degenerate, i.e. it does not lie
in any hyperplane.\\\\
$(4)$ When a projective variety $X$ is embedded in a projective
space $\P^r$ by a very ample line bundle $\mathcal{L}\in
\mbox{Pic} X $, we may write $\mathcal{O}_X (1)$ instead of
$\mathcal{L}$ so long as no confusion arise.\\\\
$(5)$ For arbitrary nonzero coherent sheaf $\mathcal{F}$ on
$\P^r$, $\mbox{Reg}(\mathcal{F})$ is defined to be
\begin{equation*}
\mbox{min}\{m|H^i (\P^r ,\mathcal{F}(m-i))=0~~~~\mbox{for
all}~~~~i\geq1\}.
\end{equation*}
Also it is well known that for an exact sequence $0 \rightarrow
\mathcal{E} \rightarrow \mathcal{F} \rightarrow \mathcal{G}
\rightarrow 0$,
\begin{enumerate}
\item[(i)] $\mbox{Reg}(\mathcal{G}) \geq
\mbox{max}\{\mbox{Reg}(\mathcal{E})-1, \mbox{Reg}(\mathcal{F})\}$
and
\item[(ii)] $\mbox{Reg}(\mathcal{E}) \geq
\mbox{max}\{\mbox{Reg}(\mathcal{F}), \mbox{Reg}(\mathcal{G})+1\}$
provided that $H^1 (\P^r ,\mathcal{E}(j))=0$ for all $j\in \Z$.
\end{enumerate}
This will be freely used without explicit mention.\\

\section{Regularity Criterion and Effects of Property $N_p$}
\noindent We begin with proving Theorem \ref{thm:regularity}.\\\\
{\bf Proof of Theorem \ref{thm:regularity}} Let $S$ be the
homogeneous coordinate ring of $\P(V)$. Since $X \hookrightarrow
\P(V)$ satisfies Property $N^S _1$, the graded $S$-module $E $
admits a minimal graded free resolution with the following
numerical type in the first two terms:
\begin{equation*}
\cdots \rightarrow S^{n_1} (-2) \rightarrow S \oplus S^t (-1)
\rightarrow E \rightarrow 0 .
\end{equation*}
Then we can prove our Theorem by using the surprising technique
used in \cite{GLP}. From the above resolution, $\mathcal{O}_X$
admits a resolution as follows:
\begin{equation*}
\cdots \rightarrow \mathcal{O}_{\P^r} ^{n_{1}} (-2)
\stackrel{v}{\rightarrow} \mathcal{O}_{\P^r} \oplus
\mathcal{O}_{\P^r} ^t (-1) \rightarrow \mathcal{O}_X \rightarrow
0.
\end{equation*}
As the proof of (Theorem 2.1, \cite{GLP}), we have the following
commutative diagram with exact rows and columns:\\
\begin{equation*}
\begin{CD}
&&&& 0          && 0 & \\
&&&& \downarrow && \downarrow & \\
&0&&& \mathcal{O}_{\P^r} &\rightarrow & \mathcal{O}_X &
\rightarrow
0\\
&\downarrow &&& \downarrow && \parallel & \\
0 \rightarrow & K & \rightarrow & \mathcal{O}_{\P^r} ^{n_1} (-2)
\stackrel{v}{\rightarrow}& \mathcal{O}_{\P^r} \oplus
\mathcal{O}_{\P^r} ^t (-1) & \rightarrow & \mathcal{O}_X &
\rightarrow
0 \\
& \downarrow && \parallel & \downarrow  && \downarrow \\
0 \rightarrow & N & \rightarrow & \mathcal{O}_{\P^r} ^{n_1} (-2)
\stackrel{u}{\rightarrow}& \mathcal{O}_{\P^r} ^t (-1) &
\rightarrow & 0 &  \\
&&&& \downarrow &&  & \\
&&&& 0          && & \\
\end{CD}
\end{equation*}
where $K$ and $N$ are kernels of $v$ and $u$, respectively. Now
let $L$  be the kernel of
\begin{equation*}
\mathcal{O}_{\P^r} \oplus \mathcal{O}_{\P^r} ^{n_0} (-1)
\rightarrow \mathcal{O}_X \rightarrow 0.
\end{equation*}
Then we can separate the above diagram as follows:\\
\begin{equation*}
\begin{CD}
              & 0             &             &  0                                                  &             &               &             \\
              & \downarrow    &             &  \downarrow                                         &             &               &             \\
0 \rightarrow & \mathcal{I}_X & \rightarrow & \mathcal{O}_{\P^r}                                  & \rightarrow & \mathcal{O}_X & \rightarrow 0\\
              & \downarrow    &             & \downarrow                                          &             & \parallel     &             \\
0 \rightarrow & L             & \rightarrow & \mathcal{O}_{\P^r} \oplus\mathcal{O}_{\P^r} ^t (-1) & \rightarrow & \mathcal{O}_X & \rightarrow 0 \\
              & \downarrow    &             & \downarrow                                          &             &               &             \\
              & \mathcal{O}_{\P^r} ^{t} (-1) &  = & \mathcal{O}_{\P^r} ^t (-1)                  & & & \\
              & \downarrow    &             & \downarrow                                          &             &               &             \\
              & 0             &             & 0                                                   &             &               &             \\
\\
\end{CD}
\end{equation*}

\begin{equation*}
\begin{CD}
              &               &             &                                 &             &   0           &             \\
              &               &             &                                 &             & \downarrow    &             \\
              & 0             &             &                                 &             & \mathcal{I}_X &             \\
              & \downarrow    &             &                                 &             & \downarrow    &             \\
0 \rightarrow & K             & \rightarrow & \mathcal{O}_{\P^r} ^{n_1} (-2)  & \rightarrow &  L            & \rightarrow 0\\
              & \downarrow    &             & \parallel                       &             & \downarrow    &             \\
0 \rightarrow & N             & \rightarrow & \mathcal{O}_{\P^r} ^{n_1} (-2)  & \rightarrow & \mathcal{O}_{\P^r} ^t (-1) & \rightarrow 0 \\
              & \downarrow    &             &                                 &             &               &             \\
              & \mathcal{I}_X &             &                                 &             &               &             \\
              & \downarrow    &             &                                 &             &               &             \\
              & 0             &             &                                 &             &               &             \\
&&&&&& \\
\end{CD}
\end{equation*}
Each diagram is directly obtained by Snake Lemma. We observe the
following facts :\\
\begin{enumerate}
\item[$(1)$] $H^1 (\P^r ,K(n))=H^1 (\P^r ,L(n))=0$ for all $n \in \Z$
by the minimality of the resolution.
\item[$(2)$] $H^2 (\P^r ,K(n))=0$ for all $n \in \Z$ by the exact
sequence
\begin{equation*}
0 \rightarrow K \rightarrow \mathcal{O}_{\P^r} ^{n_1} (-2)
\rightarrow L \rightarrow 0.
\end{equation*}
\item[$(3)$] $N$ is $(t+2)$-regular by the following Lemma
\ref{lem:ENC}.\\
\end{enumerate}
Now consider $0 \rightarrow  K \rightarrow N \rightarrow
\mathcal{I}_X \rightarrow 0$. In the exact sequence
\begin{equation*}
H^1 (\P^r , N(t+j)) \rightarrow  H^1 (\P^r , \mathcal{I}_X (t+j))
\rightarrow H^2 (\P^r , K(t+j)),
\end{equation*}
$H^1 (\P^r , N(t+j))=0$ for all $j \geq1$ and $H^2 (\P^r ,
K(t+j))=0$ for all $j\in \Z$ by the above observation. Therefore
$H^1 (\P^r ,\mathcal{I}_X (t+j))=0$ for all $j\geq1$.

For the regularity, the third condition implies that $L$ is
$(m+1)$-regular and $K$ is $(m+2)$-regular since their first
cohomology groups vanish always. So by the exact sequence $0
\rightarrow  K \rightarrow N \rightarrow \mathcal{I}_X \rightarrow
0$, we obtain that $\mathcal{I}_X $ is max$\{m+1,t+2\}$-regular.

For the third statement, consider the following commutative
diagram:
\begin{equation*}
\begin{CD}
H^0 (\P^r,N(t+2)) \otimes H^0 (\P^r,\mathcal{O}_{\P^r} (k))             & \stackrel{\alpha}{\rightarrow} &  H^0 (\P^r,N(t+2+k))   & \rightarrow 0           \\
 \gamma \downarrow                                                      &             & \downarrow   \delta        &        \\
H^0 (\P^r,\mathcal{I}_X (t+2)) \otimes H^0 (\P^r,\mathcal{O}_{\P^r}(k)) & \stackrel{\beta}{\rightarrow} &  H^0 (\P^r,\mathcal{I}_X (t+2+k)) &\\
\downarrow                                                              &             & \downarrow                  &    \\
H^1 (\P^r,K(t+2))\otimes H^0 (\P^r,\mathcal{O}_{\P^r} (k)) & & H^1
(\P^r,K(t+2+k))       &
\end{CD}
\end{equation*}
Since $N$ is $(t+2)$-regular, $\alpha$ is surjective for all $k
\geq 0$. Since $H^1 (\P^r,K(n))=0$ for all $n \in \Z$, $\gamma$
and $\delta$ are surjective. Therefore $\beta$ is surjective for
all $k \geq 0$ which completes the proof. \qed

\begin{lemma}\label{lem:ENC}
Let $\mathcal{E} = \mathcal{O}_{\P^r} ^e (-2)
\stackrel{u}{\rightarrow} \mathcal{F} = \mathcal{O}_{\P^r} ^f
(-1)$ be a generically surjective homomorphism such that
$\mbox{Supp}(\mbox{Coker}(u))$ is a finite set. Then
\begin{equation*}
Reg(Ker(u)) \leq f+2.
\end{equation*}
\end{lemma}

\begin{proof}
This can be proved by the Eagon-Northcott complex and diagram
chasing. For details, see (Lemma 5, \cite{Noma}).
\end{proof}

From now on, we proceed to prove Theorem \ref{thm:effect}. Let $X$
be a smooth complex projective variety of dimension $n$ and let
$\mathcal{L} \in \mbox{Pic}X$ be a very ample line bundle such
that
\begin{equation*}
H^1 (X,\mathcal{L}^{j})=0~~~~ \mbox{for all}~~~~j\geq2.
\end{equation*}
First we give a criterion for Property $N^S _p$.

\begin{lemma}\label{lem:criterion}
Under the hypothesis just stated, let $X \hookrightarrow \P(V)$ be
an isomorphic embedding where $V \subset H^0 (X,\mathcal{L})$ and
consider the canonical exact sequence
\begin{equation*}
0 \rightarrow \mathcal{M}_V \rightarrow V \otimes \mathcal{O}_X
\rightarrow \mathcal{L} \rightarrow 0.
\end{equation*}
For $p \geq 0$, the embedding $X \hookrightarrow \P(V)$ satisfies
Property $N^S _p$ if and only if for $0 \leq i \leq p$
\begin{enumerate}
\item[$(1)$] $\wedge^{i+1} V \otimes H^0 (X, \mathcal{L})
 \rightarrow H^0 (X, \wedge^i
\mathcal{M}_V \otimes \mathcal{L}^2 )$ is surjective and
\item[$(2)$] $H^1 (X, \wedge^{i+1} \mathcal{M}_V \otimes \mathcal{L}^j)
=0$ for all $j \geq 2$.
\end{enumerate}
\end{lemma}

\begin{proof}
This can be proved using canonical method computing Betti numbers.
Indeed assume that $E$ has a minimal free resolution of the form
\begin{equation*}
\cdots \stackrel{\varphi_{i+1}}{\rightarrow} \oplus_j
S^{k_{i,j}}(-i-j) \stackrel{\varphi_i}{\rightarrow}  \cdots
\stackrel{\varphi_1}{\rightarrow} \oplus_j S^{k_{0,j}}(-j)
\stackrel{\varphi_0}{\rightarrow} E \rightarrow 0.
\end{equation*}
As explained in \cite{Eisenbud} or \cite{Green}, $k_{i,j} =
\mbox{dim}_k \mbox{Coker}(\alpha_{i,j})$ where $\alpha_{i,j}$ is
the map contained in the exact sequence
\begin{equation*}
\begin{CD}
& \wedge^{i+1} V \otimes H^0 (X, \mathcal{L}^{j-1})
& \stackrel{\alpha_{i,j}}{\longrightarrow} & H^0 (X, \wedge^i \mathcal{M}_V \otimes \mathcal{L}^j ) \\
\rightarrow & H^1 (X, \wedge^{i+1} \mathcal{M}_V \otimes
\mathcal{L}^{j-1})& \rightarrow & \wedge^{i+1} V \otimes H^1 (X,
\mathcal{L}^{j-1}).
\end{CD}
\end{equation*}
Therefore $(1)$ is equivalent to the fact that $k_{i,2} = 0$ for
$0 \leq i \leq p$. Also for $j \geq 2$,
\begin{equation*}
k_{i,j} = h^1 (X, \wedge^{i+1} \mathcal{M}_V \otimes \mathcal{L}^
{j-1})
\end{equation*}
since $H^1 (X, \mathcal{L}^{j-1}) =0$ by assumption. So $(2)$ is
equivalent to the fact that $k_{i,j}=0$ for $0 \leq i \leq p$ and
$j \geq 2$.
\end{proof}

\begin{proposition}\label{prop:codimension1}
Let $X \subset \P(V)$ be a complex projective variety and let
$\mathcal{L}$ be a very ample line bundle on $X$. Consider
subspaces $W \subset V \subset H^0 (X,\mathcal{L})$ such that
$\mbox{codim}(W,V)=1$ and $X \rightarrow \P(W)$ is an isomorphic
embedding. Assume that $H^1 (X,\mathcal{L}^j)=0$ for all $j\geq2$.
If $X \subset \P(V)$ satisfies Property $N^S _p$ for some $p \geq
1$, then $X \subset \P(W)$ satisfies Property $N^S _{p-1}$.
\end{proposition}

\begin{proof}
Consider the following commutative diagram:\\
\begin{equation*}
\begin{CD}
              & 0             &             &   0                     &             &             &               \\
              & \downarrow    &             &  \downarrow             &             &             &               \\
0 \rightarrow & \mathcal{M}_W & \rightarrow & W \otimes \mathcal{O}_X & \rightarrow & \mathcal{L} & \rightarrow 0 \\
              & \downarrow    &             & \downarrow              &             & \parallel   &               \\
0 \rightarrow & \mathcal{M}_V & \rightarrow & V \otimes \mathcal{O}_X & \rightarrow & \mathcal{L} & \rightarrow  0 \\
              & \downarrow    &             & \downarrow              &             &             &               \\
              & \mathcal{O}_X &     =       &  \mathcal{O}_X          &             &             &               \\
              & \downarrow    &             &  \downarrow             &             &             &               \\
              & 0             &             &  0                      &             &             &
\end{CD}
\end{equation*}\\
where $\mathcal{M}_W$ and $\mathcal{M}_V$ are kernels of the above
evaluation maps. Since $X \subset \P(V)$ satisfies Property $N^S
_p$, for $0 \leq i \leq p$
\begin{enumerate}
\item[$(1)$] $\wedge^{i+1} V \otimes H^0 (X, \mathcal{L})
 \rightarrow H^0 (X, \wedge^i
\mathcal{M}_V \otimes \mathcal{L}^2 )$ is surjective and
\item[$(2)$] $H^1 (X, \wedge^{i+1} \mathcal{M}_V \otimes \mathcal{L}^j)
=0$ for all $j \geq 2$.
\end{enumerate}
by Lemma \ref{lem:criterion}. Also what we must show is that for
$0 \leq i \leq p-1$,\\
$$ (*) = \begin{cases} \wedge^{i+1} W \otimes H^0 (X, \mathcal{L})
 \rightarrow H^0 (X, \wedge^i
\mathcal{M}_W \otimes \mathcal{L}^2 )~~~~\mbox{and} \\
H^1 (X, \wedge^{i+1} \mathcal{M}_W \otimes \mathcal{L}^j)
=0~~~~\mbox{for all}~~~~ j \geq 2.
\end{cases}$$

\noindent Now consider the following commutative diagram induced
from the
above one:\\
\begin{equation*}
\begin{CD}
              & 0                          &             &        0                             &             &      0      &               \\
              & \downarrow                 &             &   \downarrow                         &             & \downarrow  &               \\
0 \rightarrow & \wedge^{i+1} \mathcal{M}_W & \rightarrow & \wedge^{i+1} \mathcal{M}_V & \rightarrow & \wedge^{i} \mathcal{M}_W  & \rightarrow 0 \\
              & \downarrow                 &             &   \downarrow                         &             & \downarrow  &               \\
0 \rightarrow & \wedge^{i+1} W \otimes \mathcal{O}_X & \rightarrow & \wedge^{i+1} V \otimes \mathcal{O}_X & \rightarrow & \wedge^{i} W \otimes \mathcal{O}_X & \rightarrow  0 \\
              & \downarrow                 &             &   \downarrow                         &             & \downarrow  &               \\
0 \rightarrow & \wedge^{i} \mathcal{M}_W  \otimes \mathcal{L} & \rightarrow & \wedge^{i} \mathcal{M}_V  \otimes \mathcal{L}   & \rightarrow &  \wedge^{i-1} \mathcal{M}_W  \otimes \mathcal{L}  & \rightarrow  0\\
              & \downarrow                 &             &   \downarrow                         &             & \downarrow  &               \\
              & 0                          &             &     0                                &             &   0         &
\end{CD}
\end{equation*}\\
This gives the following commutative diagram of cohomological groups:\\
\begin{equation*}
\begin{CD}
\wedge^{i+1} V \otimes H^0 (X,\mathcal{L}^{j-1}) & \stackrel{\gamma_{i,j}}{\rightarrow} &  \wedge^{i  } W \otimes H^0 (X,\mathcal{L}^{j-1})  & &  &   &  \\
                   @V{\alpha_{i,j}}VV                                            @VV{\beta_{i,j}}V                                                                                        &   &   \\
H^0 (X,\wedge^{i} \mathcal{M}_V  \otimes \mathcal{L}^{j}) & \stackrel{\delta_{i,j}}{\rightarrow}  & H^0 (X,\wedge^{i-1} \mathcal{M}_W  \otimes \mathcal{L}^{j})  & \rightarrow & H^1 (X,\wedge^{i } \mathcal{M}_W  \otimes \mathcal{L}^{j}) &  & \\
             \downarrow                          &                                &    \downarrow                                               &             &                                                     &             & \\
H^1 (X,\wedge^{i+1} \mathcal{M}_V  \otimes \mathcal{L}^{j-1}) &                    &  H^1 (X, \wedge^{i} \mathcal{M}_W \otimes \mathcal{L}^{j-1}) &             &                                                               &             \\
                                  &                                &    \downarrow                                                   &             &                                                     &             & \\
&                                & \wedge^{i} W \otimes H^1 (X,\mathcal{L}^{j-1})                                                 &             &                                                     &             & \\
\end{CD}
\end{equation*}\\

\noindent Note that $\alpha_{i,j}$ and $\gamma_{i,j}$ are
surjective for $j\geq2$ and $0\leq i \leq p$. So if $\delta_{i,j}$
is surjective then so is $\beta_{i,j}$. The trick is that in this
range,
\begin{equation*}
\mbox{if}~~~H^1 (X,\wedge^{i } \mathcal{M}_W  \otimes
\mathcal{L}^{j})=0,~~~~\mbox{then}~~~~H^1 (X, \wedge^{i}
\mathcal{M}_W \otimes \mathcal{L}^{j-1})=0
\end{equation*}
for all $j\geq3$. Let us remind that $H^1 (X,\mathcal{L}^\ell)=0$
for $\ell\geq 2$. This guarantees that
\begin{equation*}
H^1 (X,\wedge^{i} \mathcal{M}_W \otimes
\mathcal{L}^{j})=0~~~~\mbox{for all}~~~~j\geq
2~~~~\mbox{and}~~~~0\leq i \leq p
\end{equation*}
because $H^1 (X,\wedge^{i } \mathcal{M}_W  \otimes
\mathcal{L}^{j})=0$ for sufficiently large $j$. That is,
$\beta_{i,j}$ is surjective for all $j\geq3$ and $0\leq i \leq
p-1$. Finally $\alpha_{i,2}$ is surjective by assumption and
$\delta_{i,2}$ is surjective since it is proved that $H^1
(X,\wedge^{i} \mathcal{M}_W \otimes \mathcal{L}^{2})=0$. So
$\beta_{i,2}$ is surjective for $0\leq i \leq p-1$ which completes
the proof.
\end{proof}

\noindent {\bf Proof of Theorem \ref{thm:effect}.} For a given $V
\subset H^0 (X,\mathcal{L})$, fix a filtration
\begin{equation*}
V=V_t \subset V_{t-1} \subset \cdots \subset V_0 = H^0
(X,\mathcal{L})
\end{equation*}
by subspaces each having codimension one in the next. Then the
first statement is proved by applying Proposition
\ref{prop:codimension1} repeatedly. Also the second statement is
induced directly by Theorem \ref{thm:regularity}. \qed \\

\section{Higher Syzygies of Surface Scrolls over a Curve}
\noindent In this section we investigate Property $N_p$ of surface
scrolls over a curve and generalize Green's result about Property
$N_p$ for curves. For the convenience of the reader, we recall
notations in $\S 1$:
\begin{enumerate}
\item[$\bullet$] $C$ : smooth curve of genus $g$
\item[$\bullet$] $\mathcal{E}$ : normalized vector bundle of rank $2$ over $C$
\item[$\bullet$] $\mathfrak{e}= \wedge^2 \mathcal{E} ~~$ and $e = - \mbox{deg}(\mathfrak{e})$
\item[$\bullet$] $X = \P _C (\mathcal{E})$ : the associated ruled
surface with the projection map $\pi : X \rightarrow C $
\item[$\bullet$] $C_0$ : a minimal section such that $\mathcal{O}_X (C_0)=\mathcal{O}_{\P_C
(\mathcal{E})} (1)$
\end{enumerate}
Throughout this section we consider line bundles of the form
\begin{equation*}
\mathcal{L} =H+\pi^* B~~~~\mbox{with}~~~~\mbox{deg}(B) \geq
\mbox{max} \{ 3g,3g+e \}+1.
\end{equation*}
Note that this line bundle is very ample(Ex.5.2.11, \cite{H}) and
satisfies $H^1 (X,\mathcal{L}^j)=0$ for all $j\geq1$. Consider the
short exact sequence
\begin{equation*}
0 \rightarrow \mathcal{M}_{\mathcal{L}} \rightarrow H^0
(X,\mathcal{L})\otimes \mathcal{O}_X \rightarrow \mathcal{L}
\rightarrow 0.
\end{equation*}
Remark that $(X,L)$ satisfies Property $N_p$ if and only if $H^1
(X,\wedge^{p+1} \mathcal{M}_{\mathcal{L}} \otimes \mathcal{L})=0$
\cite{GL}.\\

\noindent {\bf Proof of Theorem \ref{thm:surfacescroll}.} Let $h^0
(X,\mathcal{L})=r+1$. For a smooth hyperplane section $C_1 (\cong
C) \in H^0 (X, \mathcal{L})$, $\mathcal{L}|_{C_1} = \mathcal{O}_C
(2\mathfrak{b}+\mathfrak{e})$ and hence we have
\begin{equation*}
0 \rightarrow  \mathcal{M}_C   \rightarrow  \mathcal{O}_C ^r
\rightarrow  \mathcal{O}_C (2\mathfrak{b}+\mathfrak{e})
\rightarrow 0.
\end{equation*}
If we pull back this sequence by $\pi$, we have the following
commutative diagram\\
\begin{equation*}
\begin{CD}
              &                       &             &                      &             & 0                              &              \\
              &                       &             &                      &             & \downarrow                     &              \\
              &         0             &             &                      &             & -C_0 +(\mathfrak{b}+\mathfrak{e})f &              \\
              &    \downarrow         &             &                      &             & \downarrow                     &              \\
0 \rightarrow & \pi ^* \mathcal{M}_C  & \rightarrow & \mathcal{O}_X ^r     & \rightarrow & (2\mathfrak{b}+\mathfrak{e})f & \rightarrow 0 \\
              & \downarrow            &             & \parallel            &             &  \downarrow                    & \\
0 \rightarrow & \mathcal{M}_{\mathcal{L}} & \rightarrow & \mathcal{O}_X ^r     & \rightarrow & \mathcal{O}_C (2\mathfrak{b}+\mathfrak{e}) & \rightarrow 0\\
              &  \downarrow           &             &                      &             & \downarrow                     &              \\
              &  -C_0 +(\mathfrak{b}+\mathfrak{e})f &            &                      &             & 0                              &              \\
              &    \downarrow         &             &                      &             &                                &              \\
              &   0                   &             &                      &             &                                &              \\
\end{CD}
\end{equation*}
by Snake Lemma where the second row is induced by the following
Lemma \ref{lem:newEuler}. From
\begin{equation*}
0 \rightarrow \pi ^* \mathcal{M}_C  \rightarrow
\mathcal{M}_{\mathcal{L}} \rightarrow \mathcal{O}_X (-C_0
+(\mathfrak{b}+\mathfrak{e})f) \rightarrow 0,
\end{equation*}
we have
\begin{equation*}
0 \rightarrow \wedge^{p+1} \pi ^* \mathcal{M}_C \rightarrow
\wedge^{p+1} \mathcal{M} \rightarrow \wedge^p \pi ^* \mathcal{M}_C
\otimes  \mathcal{O}_X (-C_0 +(\mathfrak{b}+\mathfrak{e})f)
\rightarrow 0.
\end{equation*}
and hence it is enough to check that
\begin{equation*}
H^1 (X, \wedge^{p+1} \pi ^* \mathcal{M}_C \otimes \mathcal{L})=
H^1 (X, \wedge^p \pi ^* \mathcal{M}_C \otimes
(2\mathfrak{b}+\mathfrak{e})f)= 0.
\end{equation*}
or equivalently
\begin{equation*}
H^1 (C, \wedge^{p+1} \mathcal{M}_C \otimes \mathcal{E} \otimes
\mathfrak{b})= H^1 (C, \wedge^p \mathcal{M}_C \otimes
(2\mathfrak{b}+\mathfrak{e})= 0.
\end{equation*}
Note that by using the short exact sequence $0 \rightarrow
\mathcal{O}_C \rightarrow \mathcal{E} \rightarrow \mathfrak{e}
\rightarrow 0$ corresponding to $C_0$, the first one is guaranteed
if we show $H^1 (C, \wedge^{p+1} \mathcal{M}_C \otimes B)=0$ and
$H^1 (C, \wedge^{p+1} \mathcal{M}_C \otimes \textbf{e} \otimes
B)=0$. To prove the desired vanishing we use the filtration of
$\mathcal{M}_C$. As \cite[Lemma 6]{Noma}, for generic $x_0 ,
\cdots , x_{r-2} \in C_1 \cong C$, $\mathcal{M}_C$ admits a
filtration by vector bundles
\begin{equation*}
0=\mathcal{F}^0 \subset \mathcal{F}^1 \subset \cdots \subset
\mathcal{F}^{r-1} =\mathcal{M}_C
\end{equation*}
such that
\begin{equation*}
\mathcal{F}^{r-j-1} / \mathcal{F}^{r-j-2} \cong  \mathcal{O}_C
(-B_j ) ~~~~~  (j=0, \cdots , r-2)
\end{equation*}
for effective divisors $B_j$ on $C$ whose support Supp$(B_j )$
contains $x_j$. Note that
\begin{equation*}
  B_0 + \cdots +B_{r-2}  =   2\mathfrak{b} +
\mathfrak{e}  \quad \mbox{and so} \quad \mbox{deg}(B_0) + \cdots
+\mbox{deg}(B_{r-2} ) = 2b-e.
\end{equation*}
Therefore it suffices to show that for every distinct $0 \leq i_1
, \cdots , i_{p+1} \leq r-2$,
$$\begin{cases}
H^1 (C,\mathcal{O}_C ( -B_{i_1}-\cdots - B_{i_{p+1}}
+\mathfrak{b}))=0\\
H^1 (C,\mathcal{O}_C ( -B_{i_1}-\cdots -
B_{i_{p+1}}+\mathfrak{b}+\mathfrak{e}))=0
\end{cases}$$
or equivalently, for every distinct $0 \leq j_1 , \cdots ,
j_{r-2-p} \leq r-2$,
$$\begin{cases}
H^1 (C,\mathcal{O}_C (B_{j_1}+ \cdots +B_{j_{r-2-p}}-\mathfrak{b}-\mathfrak{e}))=0\\
H^1 (C,\mathcal{O}_C (B_{j_1}+ \cdots
+B_{j_{r-2-p}}-\mathfrak{b}))=0.
\end{cases}$$
By Serre duality, we show
$$\begin{cases}
H^0 (C,\mathcal{O}_C (\mathfrak{b}+\mathfrak{e}+K_C -B_{j_1}- \cdots -B_{j_{r-2-p}}))=0\\
H^0 (C,\mathcal{O}_C (\mathfrak{b}+K_C -B_{j_1}- \cdots
-B_{j_{r-2-p}}))=0.
\end{cases}$$
Since $B_i$ is effective divisor containing $x_i$, it is enough to
prove
$$\begin{cases}
H^0 (C,\mathcal{O}_C (\mathfrak{b}+\mathfrak{e}+K_C -x_{j_1}- \cdots -x_{j_{r-2-p}}))=0\\
H^0 (C,\mathcal{O}_C (\mathfrak{b}+K_C -x_{j_1}- \cdots
-x_{j_{r-2-p}}))=0.
\end{cases}$$
This is true for generic $x_0 , \cdots , x_{r-2}$ since
$$\begin{cases}
h^0 (C,\mathcal{O}_C (\mathfrak{b}+\mathfrak{e}+K_C ))=b-e+g-1 & \leq  r-2-p=2b-e-2g-p-1\\
h^0 (C,\mathcal{O}_C (\mathfrak{b}+K_C ))=b+g-1 & \leq
r-2-p=2b-e-2g-p-1
\end{cases}$$
 which is equivalent to our
assumption
\begin{equation*}
b \geq \mbox{max}\{3g, 3g+e\}+p \quad (p \geq 0).
\end{equation*}
\qed \\

\begin{lemma}\label{lem:newEuler}
Let $Z$ be a projective variety embedded in $\P^r$ and let
$\mathcal{M} = \Omega_{\P^r} (1)|_Z $. Then for a hyperplane $H$
of $\P^r$, there is the following exact sequence:
\begin{equation*}
0 \rightarrow \mathcal{M} \rightarrow \mathcal{O}_Z ^r \rightarrow
\mathcal{O}_{Z\cap H} (1) \rightarrow 0
\end{equation*}
\end{lemma}

\begin{proof}
Let us consider the follow commutative diagram:\\
\begin{equation*}
\begin{CD}
              &             &             & 0                    &             & 0                &              \\
              &             &             & \downarrow           &             & \downarrow       &              \\
              &             &             & \mathcal{O}_Z        & =           & \mathcal{O}_Z    &              \\
              &             &             & \downarrow           &             & \downarrow       &              \\
0 \rightarrow & \mathcal{M} & \rightarrow & \mathcal{O}_Z ^{r+1} & \rightarrow & \mathcal{O}_Z (1) & \rightarrow 0 \\
              & \downarrow  &             & \downarrow           &             &  \downarrow       &               \\
0 \rightarrow & \mathcal{F} & \rightarrow & \mathcal{O}_Z ^r     & \rightarrow & \mathcal{O}_{Z \cap H} (1) & \rightarrow 0\\
              &             &             & \downarrow           &             & \downarrow       &              \\
              &             &             & 0                    &             & 0                &              \\
&&&&&&
\end{CD}
\end{equation*}
where $\mathcal{F}$ is just the kernel of $\mathcal{O}_Z ^r
\rightarrow \mathcal{O}_{Z \cap H} (1)$. By Snake lemma,
$\mathcal{F} \cong \mathcal{M}$.
\end{proof}

\section{Applications and Examples}
\subsection{Applications} Now we apply Theorem \ref{thm:effect}. First let us
recall the following works of Green\cite{Green},
Butler\cite{Butler} and Ein-Lazarsfeld\cite{EL} about Property
$N_p$:\\
\begin{enumerate}
\item[$\bullet$] (Green\cite{Green}) Let $X$ be a smooth projective
curve of genus $g$ and let $\mathcal{L}$ be a line bundle on $X$
of degree $d$. If $d\geq 2g+1+p$ for $p\geq0$, $(X,\mathcal{L})$
satisfies $N_p$.
\item[$\bullet$] (Butler,\cite{Butler}) Let $C$ be a smooth projective
curve of genus $g$. For a vector bundle $\E$ of rank $n$ over $C$,
let $X=\P(\E)$ be the associated ruled variety with tautological
line bundle $H$ and projection map $\pi:X \rightarrow C$. For a
line bundle $\mathcal{L}=aH+\pi^* B$ on $X$ with $a\geq1$ and
$B\in \mbox{PicC}$, assume that $\mu^- (\pi_* \mathcal{L})\geq
2g+2p$ for some $1\leq p\leq a$. Then $(X,\mathcal{L})$ satisfies
$N_p$.
\item[$\bullet$] (Ein-Lazarsfeld,\cite{EL}) Let $X$ be a smooth complex
projective variety of dimension $n$ with the canonical sheaf
$K_X$. For $A, B \in \mbox{Pic} X$, assume that  $A$ is very ample
and $B$ is numerically effective. Then for $p\geq0$
\begin{equation*}
K_X + (n+p)A +B ~~~~ \mbox{satisfies} ~~~~ \mbox{Property}~~~~N_p
\end{equation*}
except the case $(X,A,B)=(\P^n ,\mathcal{O}_{\P^n}
(1),\mathcal{O}_{\P^n})$ and $p=0$.\\
\end{enumerate}

\noindent Note that in these results the line bundle $\mathcal{L}$
satisfies $H^i (X, \mathcal{O}_X (j))=0$ for all $i,j\geq1$.
Therefore Theorem \ref{thm:effect}  can be applied to these
results. Recently, A. Noma obtained a very sharp regularity bound
for curves:\\
\begin{enumerate}
\item[$(*)$] (Noma,\cite {Noma})$~$  Let $X \subset \P^r (r \geq 3)$ be an
integral projective curve of degree $d$ and arithmetic genus
$\rho_a$. Let $l \leq  \mbox{min} \{ r-2 , \rho_a  \} $ be a
nonnegative integer. Then $X$ is $(d-r+2-l)$-regular and $d \geq
r+l+1$ unless $X$ is a curve embedded by a complete linear system
of degree $d
\geq 2 \rho_a +2$ and $l= \rho_a$.\\
\end{enumerate}
From this result, when $d = 2g+1+p$ for some $p \geq 1$ and $X
\subset \P^r$ by a subsystem of codimension $t \geq p-1$, then $X$
is $(t+2)$-regular and hence $(t+1)$-normal since we can take
$l=g$ and in this case $d-r+2-g=t+2$. This is the same as
information obtained by applying Theorem \ref{thm:regularity} and
Theorem \ref{thm:effect} to
Green's ``$2g+1+p$" theorem\cite{Green}. In fact, since\\
\begin{enumerate}
\item[$(**)$] (Green-Lazarsfeld,\cite{GL}) Let $\mathcal{L}$ be a line
bundle of degree $2g+p~(p\geq1)$ on a smooth projective curve $X$
of genus $g$, defining an embedding $X \hookrightarrow \P(H^0
(X,\mathcal{L}))=\P^{g+p}$. Then $\mathcal{L}$ fails to satisfy
$N_p$ if and only if either

(i) $X$ is hyperelliptic;\\
or

(ii) $X \hookrightarrow \P^{g+p}$ has a $(p+2)$-secant $p$-plane,
i.e., $H^0 (X,\mathcal{L}\otimes {K_X}^* ) \neq 0$.\\
\end{enumerate}
is known, Theorem \ref{thm:effect} refines Noma's bound a little.
We remark that in the above fact, our result give a new proof of
the relevant part:

\begin{theorem}
Let $C$ be a smooth curve of genus $g\geq 2$, and let $\mathcal{L}
\in \mbox{Pic}C$ be a line bundle of degree $d = 2g+p$. If $C
\hookrightarrow \P^{g+p}$ has a $(p+2)$-secant $p$-plane, then it
does not satisfy Property $N_{p}$.
\end{theorem}

\begin{proof}
We prove that there is a subsystem $V \subset H^0 (X,\mathcal{L})$
such that $(p+2)$-secant $p$-plane gives a $(p+2)$-secant line in
the embedding $C \hookrightarrow \P(V)=\P^{g+1}$ and hence Theorem
\ref{thm:effect} says that $N_{p}$ does not hold. Indeed
$\mbox{dim}\mbox{Sec}C =3$ and $\mbox{Sec}C \subset \P^{g+p}$ is
nondegenerate. So we can take a center of dimension $(p-1)$ in the
$p$-plane.
\end{proof}

\noindent In the same reason, Noma's bound can be refined if the
following conjecture introduced by Green and Lazarsfeld
\cite{GL} is solved:\\
\begin{enumerate}
\item[$(\star)$] Let $C$ be a smooth curve of genus $g$ and let
$\mathcal{L}$ be a very ample line bundle on $C$ with
$\mbox{deg}(\mathcal{L}) \geq 2g+1+p -2h^1
(C,\mathcal{L})-\mbox{Cliff}(C).$ Then $\mathcal{L}$ satisfies
Property $N_p$ unless $C \subset \P(H^0 (C,\mathcal{L}))$ has a
$(p+2)$-secant $p$-plane.\\
\end{enumerate}
Finally in the case of rational normal scrolls including rational
 normal curves, Property $N_p$ holds for all $p\geq0$ because
these varieties are always $2$-regular. In fact, the following is
known:\\
\begin{enumerate}
\item[$(***)$] (Ottaviani-Paoletti, \cite{OP}) The only smooth varieties in $\P^r$ such that
Property $N_p$ holds for every $p \geq 0$ are the quadrics, the
rational normal scrolls and the Veronese surface in $\P^5$.\\
\end{enumerate}
Therefore Theorem \ref{thm:effect} implies that Regularity
Conjecture is true for smooth rational scrolls:

\begin{theorem}
For a rational  scroll $X \subset \P^r$ of dimension $n$ and
degree $d$,
\begin{equation*}
\mbox{reg}(X) \leq d-(r-n)+1.
\end{equation*}
\end{theorem}

\begin{proof}
Let $\mathcal{E} = \mathcal{O}_{\P^1} (a_1 )\oplus \cdots \oplus
\mathcal{O}_{\P^1} (a_n )$ be such that $a_i \geq1$ and let $X =
\P_{\P^1} (\mathcal{E})$ with the tautological line bundle
$\mathcal{L}$. Recall that the degree $d$ of $\mathcal{L}$ is
equal to $a_1 +\cdots +a_n$, $\mathcal{O}_X$ is $3$-regular, and
$h^0 (X,\mathcal{L})=a_1 +\cdots +a_n +n$. When $X \subset \P^r$
by $\mathcal{L}$, let $t=h^0 (X,\mathcal{L})-r-1$. Then $X$
satisfies max$\{3,t+2\}$-regular. But $t+2=d-(r-n)+1$. Also if
$t=0$, then clearly $X$ is $2$-regular since the genus of $\P^1$
is $0$.
\end{proof}

\noindent Remark that Bertin\cite{Bertin} proved Regularity
Conjecture for smooth ruled scrolls over a curve of arbitrary
genus using different method.

\subsection{Examples}

Here we give concrete examples which say that the range of $t$
cannot be refined. As these examples show, we can also know that
the failure of higher order normality of noncomplete embedding
gives some obstruction to Property $N_p$ for linearly normal
embedding.\\

\noindent {\bf Example 1.} Let $C$ be a smooth projective curve of
genus $2$ and $B \in \mbox{Pic}C$ a line bundle of degree $6$. Now
consider $C \subset \P^3$ given by a subsystem of $H^0 (C, B)$.
Then  $C$ does not satisfy $2$-normality since $h^0 (\P^3
,\mathcal{O}_{\P^3} (2)) = 10 < h^0 (C , 2B)=11$. So it fails to
be $3$-regular. In fact it is exactly $4$-regular by (Theorem 1,
\cite{Noma}).  Note that $B$ satisfies Property $N_1$ by
\cite{Green} and fails to satisfy $N_2$ by Theorem
\ref{thm:effect} or by \cite{GL} since $C$ is hyperelliptic.\\

\noindent {\bf Example 2.} For $C$ and $B$ in Example 1, let
$\mathcal{E} = \mathcal{O}_C (B) \oplus \mathcal{O}_C (B) $ and $X
= \P_C (\mathcal{E})$ with the tautological line bundle
$\mathcal{O}_X (1)$. Observe the followings:
\begin{enumerate}
\item[(1)]
For the Segre embedding $Y = \P^4 \times \P^1 \hookrightarrow
\P^9$, choose $P=(P_1 , P_2 ), Q=(Q_1 , Q_2 ) \in Y$ such that $
P_i \neq Q_i $ for $i=1,2$. Then $\overline{PQ} \cap Y =\{P, Q\}$.
\item[(2)] Since $X \cong C \times \P^1$, $C \stackrel{\mid B
\mid}{\rightarrow} \P^4$ gives $X \hookrightarrow Y
\hookrightarrow \P^9$ which is just the embedding by
$\mathcal{O}_X (1)$. In this case, $\mbox{Sec} (X) \neq Y$ by
$(1)$. So we can choose a point $\alpha \in Y - Sec (X)$. Let $f :
X \rightarrow \P^8$ be the linear projection from $\alpha$.
\end{enumerate}
Now we claim that $f(X) \subseteq \P^8$ fails to be $3$-regular.
Let $\Delta \subset \P^9$ be the fibre of $\P^4 \times \P^1
\rightarrow \P^1$ containing $\alpha$. Then $f$ maps $\Delta$ to a
$3$-dimensional linear subspace  $\Lambda \subset \P^8$. Let
$f(X\cap \Delta)$ be $C_1 \cong C$. Then we have the following
commutative diagram:
\begin{equation*}
\begin{CD}
& 0 &  & 0 & & 0 & \\
& \downarrow &  & \downarrow & & \downarrow & \\
0 \rightarrow & \mathcal{I}_X \cap \mathcal{I}_\Lambda &
\rightarrow & \mathcal{I}_\Lambda &  \rightarrow & \mathcal{O}_X
(-C_1 ) & \rightarrow 0 \\
& \downarrow &  & \downarrow & & \downarrow & \\
0 \rightarrow& \mathcal{I}_X &\rightarrow& \mathcal{O}_{\P^8} &
\rightarrow& \mathcal{O}_X & \rightarrow 0 \\
& \downarrow &  & \downarrow & & \downarrow & \\
0 \rightarrow & \mathcal{I}_{C_1} & \rightarrow &
\mathcal{O}_{\P^3} &
\rightarrow & \mathcal{O}_{C_1} & \rightarrow 0\\
& \downarrow &  & \downarrow & & \downarrow & \\
& 0 &  & 0 & & 0 &\\
&&&&&&
\end{CD}
\end{equation*}
where $\mathcal{I}_X ,\mathcal{I}_{\Lambda}$ and
$\mathcal{I}_{C_1}$ are the ideal sheaves of  $X, \Lambda \subset
\P^8$ and $C_1 \subset \Lambda \cong \P^3$, respectively. Note
that in general $ \mathcal{I}_X \rightarrow \mathcal{I}_{C_1}$ is
not surjective. But in our case, this is surjective since $X$ is a
product of other varieties and the embedding is just Segre
embedding. From this diagram, if $h^1 (\P^8 ,\mathcal{I}_X (2)
)=0$, then $h^1 (\P^3 ,\mathcal{I}_{C_1} (2) )=0$ since $h^2
(\Lambda,\mathcal{I}_X \cap \mathcal{I}_\Lambda (2))=0$. But as
explained in Example 1,  $C_1 \subset \P^3$ can never be
$2$-normal. Therefore $h^1 (\P^8 ,\mathcal{I}_X (2) )\neq 0$ and
the regularity of $X \subseteq \P^8$ is at least $4$. Note that
since $\mu^- (\mathcal{E})=6$, $\mathcal{O}_X (1)$ has Property
$N_1$ by \cite{Butler}. On the other hand Theorem \ref{thm:effect}
implies the failure of Property $N_2$.\\

\noindent {\bf Example 3.} It is well known and can be seen in
\cite{OP} that
$$(\P^2 ,\mathcal{O}_{\P^2} (d)) ~~~~ \mbox{has Property $N_p$}~~~~ \begin{cases} ~~~~\mbox{for all}~~~~p\geq0~~~~\mbox{when}~~~~d=2,~~~~\mbox{and} \\
                                 \mbox{if and only if}~~~~p \leq 3d-3 ~~~~\mbox{when}~~~~d\geq3.
\end{cases}$$
In particular for $(\P^2 ,\mathcal{O}_{\P^2} (3))$, our results
about Property $N^{S}_p$ can be applied to every subsystem giving
embedding. For instance, subsystems of $H^0 (\P^2
,\mathcal{O}_{\P^2} (3))$ can be understood as the following
table. Note that this is optimal since there exists a center such
that $\P^2 \rightarrow \P^{9-t}$ has a $(t+1)$-secant line.
\begin{table}[hbt]
\begin{center}
\begin{tabular}{|c|c|c|c| }\hline
   $t$      &  Embedding               & Higher normality                & Regularity           \\ \hline \hline
   $0$      &  $\P^2 \rightarrow \P^9$ &  projectively normal            &  $2$-regular         \\ \hline
   $1$      &  $\P^2 \rightarrow \P^8$ & $k$-normal for all $k \geq 2$   &  $3$-regular         \\ \hline
   $2$      &  $\P^2 \rightarrow \P^7$ & $k$-normal for all $k \geq 3$   &  $4$-regular         \\ \hline
   $3$      &  $\P^2 \rightarrow \P^6$ & $k$-normal for all $k \geq 4$   &  $5$-regular         \\ \hline
   $4$      &  $\P^2 \rightarrow \P^5$ & $k$-normal for all $k \geq 5$   &  $6$-regular         \\ \hline
\end{tabular}
\end{center}
\caption{The case of $(\P^2 ,\mathcal{O}_{\P^2} (3))$}
\end{table} \\

\noindent {\bf Remark 1.} For a smooth projective variety $X$,
assume that $H^i (X,\mathcal{L}^{j})=0$ for all $i,j \geq 1$ and
$H^2 (X,\mathcal{O}_X )=H^n (X,\mathcal{O}_X (2-n))=0$(e.g., $X$
is a curve or a ruled scroll over a curve etc). Then
$\mathcal{O}_X$ is $2$-regular thanks to the Kodaira vanishing
theorem. Property $N_p$ for those varieties implies that $X\subset
\P(V_t )$ is max$\{3,t+2\}$-regular for $0 \leq t \leq p-1$.\\

\noindent {\bf Remark 2.} It seems very interesting to find a
range of $t$ such that Property $N^{S}_1$ holds since this
guarantees that the regularity or higher order normality changes
as expected under linear projections without depending on the
center. For ruled varieties, it looks possible to investigate this
kind of expected good behaviors by analyzing the vector bundle
$\mathcal{M}_V$ \cite{KP}. This leads to the connection between
geometry of higher secant varieties and higher linear syzygies of
a smooth projective variety in some sense.


\begin{thebibliography}{00}
\bibitem{AR} Alberto Alzati and Francesco Russo, {\em On the
$k$-normality of projected algebraic varieties}, Bull. Braz. Math.
Soc. (N.S.) 33 (2002), no. 1, 27--48.

\bibitem{BM} Dave Bayer and David Mumford,
{\em What can be computed in algebraic geometry?}, Computational
algebraic geometry and commutative algebra, Cambridge Univ. Press
Cambridge, (1993), pp. 1-48.

\bibitem{Bertin} Marie-Am\'{e}lie Bertin,
{\em On the regularity of varieties having an extremal secant
line}, J.reine angew.Math. 545 (2002), 167-181.

\bibitem{Birkenhake} Christina Birkenhake, {\em Noncomplete linear systems on abelian varieties},
Trans. Amer. Math. Soc. 348 (1996), no. 5, 1885--1908.

\bibitem{Butler} David C.Butler,
{\em Normal generation of vector bundles over a curve},
J.Differential Geo,39 (1994), 1-34.

\bibitem{C} G. Castelnuovo,
{\em Sui multipli di une serie lineare di gruppi di punti
appartenente ad une curva algebraic}, Rend. Circ. Mat. Palermo (2)
7 (1893), 89-110.

\bibitem{E} D. Eisenbud,
{\em Commutative Algebra with a view Toward Algebraic Geometry},
no. 150, Springer-Velag New York, (1995)

\bibitem{Eisenbud} Daivid Eisenbud,
{\em The Geometry of Syzygies}, Lecture Note, (2001)

\bibitem{EG} D. Eisenbud and S. Goto,
{\em Linear free resolutions and minimal multiplicity}, J. Alg. 88
(1984), 89-133.

\bibitem{EL} Lawrence Ein and R. Lazarsfeld,
{\em Syzygies and Koszul cohomology of smooth projective varieties
of arbitrary dimension}, Inv. Math. 111 (1993), 51-67.

\bibitem{Green} M.Green,
{\em Koszul cohomology and the geometry of projective varieties},
J.Differential Geo,19 (1984), 125-171.

\bibitem{GL} M.Green and R. Lazarsfeld,
{\em Some results on the syzygies of finite sets and algebraic
curves}, Composition Mathematica,67 (1988), 301-314.

\bibitem{GLP} L. Gruson, R. Lazarsfeld and C. Peskine,
{\em On a theorem of Castelnuovo and the equations defining
projective varieties}, Inv. Math. 72 (1983), 491-506.

\bibitem{H} Robin Hartshorne,
{\em Algebraic Geometry}, no. 52, Springer-Velag New York, (1977)

\bibitem{K} Sijong Kwak,
{\em Castelnuovo Regularity of Smooth Projective Varieties of
Dimension 3 and 4},  J. Alg. Geom. 7(1998), 195-206.

\bibitem{KP} Sijong Kwak and Euisung Park, {\em Defining
equations and higher normality of ruled varieties over a curve},
preprint.

\bibitem{L} R. Lazarsfeld,
{\em A sharp Castelnuovo bound for smooth surfaces}, Duke Math. J.
55 (1987), no. 2, 423--429.

\bibitem{M} David Mumford,
{\em Lectures on curves on an algebraic surface},  Ann. of Math.
Studies. 59(1966).

\bibitem{M2} David Mumford,
{\em Varieties defined by quadric equations},  Corso CIME in
Questions on Algebraic Varieties, Rome, 30-100(1970)

\bibitem{Noma} Atsushi Noma,
{\em A bound on the Castelnuovo-Mumford regularity for curves},
Math.Ann.322, 69-74 (2002)

\bibitem{OP} G.Ottaviani and R.Paoletti {\em Syzygies of Veronese
embeddings}, Composition Mathematica,125 (2001), 31-37.

\end{thebibliography}
\end{document}